\documentclass[11pt]{article} 
\RequirePackage{amsmath}
\RequirePackage{amsfonts}
\usepackage[latin1]{inputenc}
\usepackage[T1]{fontenc}
\usepackage[french]{babel}

\makeatletter
\def\@sect#1#2#3#4#5#6[#7]#8{%
  \ifnum #2>\c@secnumdepth
    \let\@svsec\@empty
  \else
    \refstepcounter{#1}%
    \protected@edef\@svsec{\@seccntformat{#1}\relax}%
  \fi
  \@tempskipa #5\relax
  \ifdim \@tempskipa>\z@
    \begingroup
      #6{%
        \@hangfrom{\hskip #3\relax\@svsec}%
          \interlinepenalty \@M #8\@@par}%
    \endgroup
    \csname #1mark\endcsname{#7}%
    \addcontentsline{toc}{#1}{%
      \ifnum #2>\c@secnumdepth \else
        \protect\numberline{\csname the#1\endcsname.}%
      \fi
      #7}%
  \else
    \def\@svsechd{%
      #6{\hskip #3\relax
      \@svsec #8}%
      \csname #1mark\endcsname{#7}%
      \addcontentsline{toc}{#1}{%
        \ifnum #2>\c@secnumdepth \else
          \protect\numberline{\csname the#1\endcsname.}%
        \fi
        #7}}%
  \fi
  \@xsect{#5}}
\def\@seccntformat#1{\csname the#1\endcsname.\quad}

\def\@begintheorem#1#2{\trivlist
   \item[\hskip \labelsep{\bfseries #1\ #2.}]\itshape}
\def\@opargbegintheorem#1#2#3{\trivlist
      \item[\hskip \labelsep{\bfseries #1\ #2\ (#3).}]\itshape}

\renewcommand\theequation{\thesection.\arabic{equation}}
\@addtoreset{equation}{section}
\makeatother

\newtheorem{theo}[equation]{Th\'eor\`eme}
\newtheorem{cor}[equation]{Corollaire}
\newtheorem{question}[equation]{Question}
\newenvironment{exemple}{
\refstepcounter{equation}\trivlist%
\item[\hskip \labelsep{\bfseries Exemple \theequation.\ }]}%
{\endtrivlist}%
\newenvironment{remarque}{
\refstepcounter{equation}\trivlist%
\item[\hskip \labelsep{\bfseries Remarque \theequation.\ }]}%
{\endtrivlist}%

\newcommand{\pour}{\mathrm{\ pour\ }}
\newcommand{\de}{\mathrm{ d }}
\newcommand{\oper}[2]{\newcommand{#1}{\mathop{\mathrm{#2}}\nolimits} }
\oper{\Ric}{Ric}
\oper{\diam}{diam}
\oper{\injrad}{inj}
\oper{\Vol}{Vol}
\oper{\SL}{SL}
\oper{\SO}{SO}
\oper{\dimension}{dim}
\newcommand{\N}{\mathbb N}
\newcommand{\Z}{\mathbb Z}
\newcommand{\R}{\mathbb R}

\DeclareSymbolFont{greek}{OML}{ptmcm}{m}{it}
\DeclareMathSymbol{\codiff}{\mathord}{greek}{"0E}

\title{Effondrements et petites valeurs
propres des formes différentielles}
\author{Pierre Jammes}
\date{}
\begin{document}
\maketitle
{\small 
\textsc{Résumé.---}
À courbure et diamètre bornés, les valeurs propres non nulles du 
laplacien de Hodge-de~Rham agissant sur les formes différentielles 
d'une variété compacte ne sont pas uniformément minorées comme c'est 
le cas pour les fonctions, et si l'une d'elle tend vers zéro alors
le volume de la variété tend aussi vers zéro, c'est-à-dire qu'elle
s'effondre. On présente ici les résultats obtenus ces dernières années
concernant le problème réciproque, à savoir déterminer le comportement
asymptotique des premières valeurs propres d'une variété lorsqu'elle 
s'effondre. 

Mots-clefs : effondrement de variétés, laplacien, formes différentielles, 
petites valeurs propres.

MSC2000 : 58J50, 58C40}

\section{Introduction}
Soit $(M^n,g)$ une vari\'et\'e riemannienne compacte connexe orientable
de dimension $n$. Le laplacien de Hodge-de~Rham, agissant sur l'espace
$\Omega^p(M)$ des $p$-formes différentielles de $M$, est défini par
$\Delta=\de\codiff+\codiff\de$, où $\de$ désigne la différentielle extérieure
et $\codiff$ la codifférentielle, adjoint de $\de$ pour le
produit scalaire $L^2$ sur $M$. Lorsque $p=0$, on retrouve le 
laplacien agissant sur les fonctions.

Le spectre du laplacien de Hodge-de~Rham forme un ensemble discret de 
nombres positifs ou nuls qu'on notera
$$0=\lambda_{p,0}(M,g)<\lambda_{p,1}(M,g)\leq\lambda_{p,2}(M,g)\leq\dots,$$
où les valeurs propres non nulles sont répétées s'il y a multiplicité.
La multiplicité de la valeur propre nulle est un invariant topologique:
c'est en fait le $p$-ième nombre de Betti de la variété $M$.

L'étude du spectre du laplacien agissant sur les fonctions montre
qu'à diamètre borné et courbure de Ricci minorée, la 1\iere{} valeur
propre ne peut pas être arbitrairement petite :
\begin{theo}\label{intro:th1}
Soit $a$ et $d$ deux réels strictement positifs et $n$ un entier.
Il existe une constante $c(n,a,d)>0$ telle que si $(M,g)$ est
une variété riemannienne de dimension $n$ dont le diamètre et la courbure
de Ricci vérifient $\diam(M,g)\leq d$ et $\Ric(M,g)\geq -ag$, alors
$$\lambda_{0,1}(M,g)\geq c.$$
\end{theo}
Notons que ce problème a fait l'objet de nombreux travaux (\cite{ly80}, 
\cite{gr80}, \cite{bbg85} et \cite{cz95}, par exemple).

B.~Colbois et G.~Courtois ont montré dans \cite{cc90} que le théorème 
\ref{intro:th1}
ne se généralise pas aux formes différentielles pour $0<p<n$
---~même avec l'hypothèse plus forte de courbure sectionnelle bornée~---
en donnant pour tout $p$ des exemples de variétés $M$ admettant
une suite de métriques $(g_i)$ telle que 
$\lim_{i\to\infty}\lambda_{p,1}(M,g_i)=0$, la courbure et le diamètre
de $M$ étant uniformément bornés. En outre, ils mettent en évidence
le fait que si une valeur propre tend vers zéro à courbure et diamètre
bornées ---~ce qu'on appelera \emph{petite valeur propre}~---
alors le volume de la variété (ou de manière équivalente son rayon 
d'injectivité) tend aussi vers zéro, c'est-à-dire qu'elle
s'effondre. Ces résultats soulèvent deux questions :
\begin{question}\label{intro:q1}
À quelles conditions une variété qui s'effondre
admet-elle une ---~ou plusieurs~--- petite valeur propre?
\end{question}
\begin{question}\label{intro:q2}
À quelle vitesse les petites valeurs propres tendent-elles vers zéro
par rapport au volume ou au rayon d'injectivité de la variété?
\end{question}

\section{Contexte topologique et géométrique}
Toutes les variétés n'admettent pas d'effondrement à courbure bornée,
il existe des obstructions topologiques. On peut par exemple
remarquer que, par définition, une variété qui s'effondre a un
volume minimal nul, et par conséquent sa caractéristique d'Euler
et tous ses nombres caractéristiques sont nuls.

L'étude de ces effondrement a fait l'objet de nombreux travaux, dont on peut
par exemple trouver un présentation sythétique dans \cite{ro02}, et on 
peut en fait décrire précisément les variétés qui s'effondrent.
Si on se donne une suite de métrique qui effondre une variété $M$,
on peut en extraire une sous-suite $(g_i)$ telle que $(M,g_i)$
converge pour la distance de Gromov-Hausdorff vers un espace
métrique.  Dans le cas où cet espace  est une
variété riemannienne lisse $(N,h)$, on sait que $M$ possède une structure 
de fibré sur $N$ :
\begin{theo}[\cite{fu87}]\label{topo:th}
Soit $(M_i,g_i)$ une suite de variétés compactes de dimension $n$ et
$(N,h)$ une variété riemannienne compacte de dimension $m<n$. Si pour 
tout $i$ la courbure sectionnelle de $M_i$ vérifie $|K(M,g_i)|\leq 1$,
et si $(M,g_i)$ converge 
vers $(N,h)$ pour la distance de Gromov-Hausdorff, alors pour tout
$i$ suffisamment grand il existe une fibration $\pi_i:M_i\rightarrow N$ 
dont la fibre est une infranilvariété.
\end{theo}
K.~Fukaya montre en outre dans \cite{fu89} un résultat semblable dans le cas 
général, l'espace limite de l'effondrement étant alors une variété stratifiée.

 Dans \cite{cfg92}, J.~Cheeger, K.~Kukaya et M.~Gromov étudient
de manière précise la métriques des variétés effondrées et montrent
que dans la situation du théorème \ref{topo:th}, la métrique le
long des fibres est proche d'une métrique invariante pour la
structure nilpotente de la fibre.

 Pour l'étude des petites valeurs propres, on se place donc souvent
dans une situation  $M$ possède une structure de fibré, en
se restreignant eventuellement au cas plus simple où la base ou la 
fibre est fixée.

\section{Cohomologie limite}
Pour répondre à la question \ref{intro:q1} dans le cas où $(M,g_i)$
tend vers une variété riemannienne $(N,h)$, J.~Lott définit dans \cite{lo02} 
un opérateur limite pour le laplacien: quitte à extraire une 
sous-suite de $(g_i)$, le spectre du laplacien de Hodge-de~Rham converge 
pour tout $p$ et 
il existe un opérateur noté $\Delta_\infty^p$ agissant sur un espace des 
formes différentielles $\Omega^p(N,E^*)$, où $E^*$ est un fibré vectoriel
gradué 
sur $N$ dépendant de $M$, et dont le spectre $\{\lambda_{p,k}^\infty\}$
vérifie $\lambda_{p,k}^\infty=\displaystyle
\lim_{i\to+\infty}\lambda_{p,k}(M,g_i)$.
La multiplicité de la valeur propre nulle de $\Delta_\infty^p$ donne
donc le nombre de valeurs propres petites ou nulles produites par
l'effondrement. Par théorie de Hodge, on peut aussi identifier 
le noyau de ce laplacien limite à une cohomologie limite de $\Omega^*(N,E^*)$.

 L'étude de ces objets limites n'est pas aisée, et les résultats
généraux de \cite{lo02} étant très techniques, nous ne les énoncerons 
pas ici. Nous allons cependant en donner des corollaires simples.

 Une situation assez bien comprise est celle des fibrés en tores sur le 
cercle s'effondrant sur leur base. Elle illustre le rôle de la topologie
dans l'existence de petites valeurs propres. En combinant les arguments
de \cite{lo02} et \cite{jath} (voir aussi \cite{ja03}), on peut écrire :

\begin{theo}\label{coh:th2}
Soit $n\geq2$, $A\in\SL_n(\Z)$, $d$ la multiplicité algébrique de la 
valeur propre $1$ de $A$ et $d'$ sa multiplicité géométrique.
On considère le fibré $M\stackrel\pi\rightarrow S^1$ de fibre $T^n$ construit
par suspension du difféomorphisme $A$. 
 Alors:
\begin{enumerate}
\item $b_1(M)=d'+1$; 
\item
Si $(g_i)_{i\in\N}$ est une suite de métrique sur $M$ telle que
$|K(M,g_i)|\leq a$ et que la submersion $\pi$ soit une
$\frac1i$-approximation de Hausdroff, alors il existe une constante
$c(M,a)>0$ telle que $\lambda_{1,d-d'+1}(M,g_i)>c$ pour tout $i$.
\item
Pour tout $k\leq d-d'$, il existe une famille de  métriques 
$(g^k_\varepsilon)$ sur $M$ de courbure et diamètre uniformément 
bornés par rapport à $\varepsilon$ et une constante $c''(M)>0$
telle que $(M,g_\varepsilon)\stackrel{d_\textrm{G-H}}{\longrightarrow}S^1$ et
$\lambda_{1,i}(M,g^k_\varepsilon)\rightarrow 0$ pour 
$i\leq k$ quand $\varepsilon\rightarrow 0$, et
$\lambda_{1,k+1}(M,g^k_\varepsilon)>c''$ si $k<n$.
\item Si la matrice $A$ est semi-simple, alors $M$ n'admet pas de
petites valeurs propres en s'effondrant sur $S^1$, quel que soit le
degré.
\end{enumerate}
\end{theo}

\begin{remarque}
Les points 2 et 4 mettent en évidence le rôle de la topologie dans 
l'existence de petites valeurs propres, et le point 3 montre que
la géométrie de l'effondrement a elle aussi une influence.
\end{remarque}

Pour les effondrements sur une variété $N$ quelconque, un autre résultat de 
\cite{lo02} donne des majorations du nombre de 
petites valeurs pour les $1$-formes en fonction de la topologie de $M$ et $N$:
\begin{theo}[\cite{lo02}]
Lorsque $M$ s'effondre à courbure bornée sur $N$, le nombre
maximal $m$ de petites valeurs propres non nulles pour les $1$-formes
vérifie
\begin{equation}
m\leq \dimension(M)-\dimension(N)+b_1(N)-b_1(M)
\end{equation}
et
\begin{equation}\label{coh:eq2}
m\leq\dimension(M).
\end{equation}
\end{theo}
L'inégalité (\ref{coh:eq2}) a la particularité de ne pas dépendre
de l'espace limite, mais elle est peu précise
car à dimension fixée, la topologie influe sur le nombre de petites valeurs 
propres comme le montre le théorème \ref{coh:th2}, ce qui soulève le
problème suivant :
\begin{question}
Comment estimer le nombre maximal de petites valeurs propres que
peut admettre une variété $M$ donnée ---~indépendamment de l'espace
limite de l'effondrement~--- en fonction de sa topologie?
\end{question}

\section{Minoration du spectre}
\subsection{Minoration par le rayon d'injectivité}
Selon \cite{cc90}, si le diamètre, la courbure et le rayon d'injectivité 
d'une variété compacte $(M^n,g)$ vérifient $|K(M,g)|<a$, $\diam(M,g)<d$ et
$\injrad(M,g)>r$, où $a$, $d$ et $r$ sont des réels strictement positifs,
alors il existe une constante $c(n,a,d,r)>0$ telle que 
$\lambda_{p,1}(M,g)>c(n,a,d,r)>0$, pour tout $p$. La démonstration
donnée dans \cite{cc90} ne permet malheureusement pas d'expliciter cette
constante. Un progrès notable a été réalisé par S.~Chanillo et F.~Trèves
dans \cite{ct97} en donnant une minoration explicite du spectre en fonction
du rayon d'injectivité et du nombre de boules géodésiques permettant
de recouvrir $M$. Plus précisément, on se donne un réel $0<r<\injrad(M,g)$,
et on note $N$ le nombre de boules géodésiques de rayon $4^{-n}r$ 
nécessaires pour recouvrir la variété $M$. On a alors:
\begin{theo}[\cite{ct97}] Il existe une constante $C(n,a)>0$ telle que 
si $|K(M,g)|<a$, alors
$$\lambda_{p,1}(M,g)\geq C\cdot r^{-2}N^{-4(n+1)},$$
pour tout $p$.
\end{theo}
Comme l'ont remarqué B.~Colbois et G.~Courtois dans \cite{cc00}, on peut 
en déduire une minoration en fonction du seul rayon d'injectivité :
\begin{cor}\label{min:cor1}
Pour tous réels $a$ et $d$ strictement positifs, il existe une constante 
$c(n,a,d)>0$ telle que si $|K(M,g)|<a$ et $\diam(M,g)<d$, alors
$$\lambda_{p,1}(M,g)\geq c\cdot\injrad(M,g)^{4n^2+4n-2},$$
pour tout $p$.
\end{cor}

\begin{remarque}
Ces estimations sont très générales: il n'y a aucune hypothèse sur
la topologie de la variété.
\end{remarque}

Nous allons démontrer ici qu'un argument élémentaire de théorie de Hodge 
permet de préciser ces minorations, en particulier lorsque $p$ est
petit par rapport à $n$.
\begin{theo}\label{min:th2}
Pour tous réels $a$ et $d$ strictement positifs,
Il existe des constantes $C_1(n,a),C_2(n,a,d)>0$ telle que
si $|K(M,g)|<a$ et $\diam(M,g)<d$, alors on a pour tout $p$
$$\lambda_{p,1}(M,g)\geq C_1\cdot r^{-2}N^{-7(p+1)}$$
et
$$\lambda_{p,1}(M,g)\geq C_2\cdot\injrad(M,g)^{7n(p+1)-2}.$$
\end{theo}
\textbf{Démonstration:}
Soit $\lambda$ une valeur propre non nulle du laplacien agissant sur les 
$p$-formes.
Si $\lambda$ admet une forme propre exacte $\omega$, S.~Chanillo et F.~Trèves
exhibent dans la démonstration du
théorème 1.1 de \cite{ct97} une $(p-1)$-forme $\varphi$ et une constante
$c(n,a)>0$ telle que $\de\varphi=\omega$ et $\|\varphi\|_2\leq c
\cdot rN^{\frac72p}\|\omega\|_2$. Si on note $\varphi'$ la forme
propre coexacte telle que $\de\varphi'=\omega$, on a alors  
$\lambda=\frac{\|\omega\|_2^2}{\|\varphi'\|_2^2}
\geq\frac{\|\omega\|_2^2}{\|\varphi\|_2^2}\geq 
c^{-2}r^{-2}N^{-7p}.$
Si $\lambda$ n'admet pas de $p$-forme propre exacte, alors elle
admet une $p$-forme propre coexacte, dont la différentielle sera
une $(p+1)$-forme propre exacte, de même valeur propre. En appliquant 
l'inégalité précédente à cette $(p+1)$-forme, on obtient, $\lambda\geq
c^{-2}r^{-2}N^{-7(p+1)}$.

Selon \cite{gr80}, le nombre $N$ est majoré par 
$c'\injrad(M,g)^{-n}$, où $c'$ est
une constante strictement positive dépendant de $a$, $d$ et $n$, ce
qui permet de conclure à la minoration 
$\lambda_{p,1}(M,g)\geq C_2\cdot\injrad(M,g)^{7n(p+1)-2}$.
\hfill$\square$

\begin{remarque}
Quand $p>\frac n2$, on peut améliorer la minoration en utilisant 
le fait que $\lambda_{p,1}(M,g)=\lambda_{n-p,1}(M,g)$.
\end{remarque}

À degré fixé, l'exposant du rayon d'injectivité dans la minoration
du théorème \ref{min:th2}
est une fonction affine de $n$ au lieu d'une fonction quadratique comme
au corollaire \ref{min:cor1}. On peut envisager de se débarrasser de cette 
dépendance par rapport au degré:
\begin{question}\label{min:q1}
Peut-on obtenir une minoration du spectre de la forme $\lambda_{p,1}(M,g)>
c(n,a,d)\cdot\injrad(M,g)^{an+b}$, où $a$ et $b$ sont des réels
indépendants de $p$ ?
\end{question}

\subsection{Spectre des fibrés principaux en tores}
 Une situation pour laquelle le comportement asymptotique des petites
valeurs propres est bien connu est celle des fibrés en cercles
s'effondrant sur leur base, qui a été étudiée par B.~Colbois et G.~Courtois:
\begin{theo}[\cite{cc00}]\label{min:cercle}
Soit $a$ et $d$ deux réels strictement positifs et $M\stackrel\pi\to
N$ un fibré en cercle de dimension $n$ et de classe d'Euler $[e]\in H^2(N)$ 
sur une variété riemannienne $(N,h)$. Il existe
des constantes $\varepsilon_0(n,a,d,(N,h))>0$ et $C_i(n,a,d,(N,h))>0$
pour $i=1,2,3$ telles que si $g$ est une métrique sur $M$ 
vérifiant $\diam(M,g)\leq d$, $|K(M,g)|\leq a$ et telle que la submersion
$\pi$ soit une $\varepsilon$-approximation de Hausdorff pour 
$\varepsilon\leq\varepsilon_0$, alors pour
$1\leq p\leq n$,
\begin{enumerate}
\item\label{vol:th1:1}$\displaystyle\lambda_{p,m_p+1}(M,g)\geq C_1$;
\item Si $[e]\neq0$, alors $\displaystyle C_2\|e\|_2^2\varepsilon^2\leq
\lambda_{1,1}(M,g) \leq C_3\|e\|_2^2\varepsilon^2$, où $\|e\|_2$
est la norme du représentant harmonique de $[e]$;
\item\label{vol:th1:3} Si $\dim H^2(N,\R)=1$, alors
$$ C_2\|e\|_2^2\varepsilon^2\leq
\lambda_{p,k}(M,g) \leq C_3\|e\|_2^2\varepsilon^2 \pour 1\leq k\leq m_p,$$
\end{enumerate}
avec $m_p=b_p(N)+b_{p-1}(N)-b_p(M)$.
\end{theo}
 On peut noter que l'exposant de $\varepsilon$ dans les estimations
des petites valeurs propres est indépendant de la dimension, et que le 
paramètre $\varepsilon$ est l'ordre de grandeur de la longueur 
des fibres, qu'on peut interpréter comme l'ordre de grandeur du 
rayon d'injectivité ou du volume de $M$. La question se pose
naturellement de savoir dans quelle mesure ces estimations se généralisent:
\begin{question}\label{min:q2}
Peut-on obtenir une minoration de $\lambda_{p,1}(M,g)$ asymptotiquement 
de l'ordre de $\injrad(M,g)^2$ ou $\Vol(M,g)^2$ quand la variété s'effondre ?
\end{question}
Notons qu'une minoration du spectre par $\Vol(M)^2$ entrainerait
une minoration par $\injrad(M)^{2n}$, ce qui répondrait affirmativement
à la question \ref{min:q1}. Par ailleurs, on ne peut pas en général
majorer les premières valeurs propres comme dans le théorème \ref{min:cercle}:
quand la fibre est de dimension strictement supérieure à 1, le
nombre de petites valeurs propres varie \emph{a priori} avec 
la géométrie de l'effondrement,
même en supposant que le fibré est principal (voir \cite{ja03}).

 L'étude des fibrés principaux en tores menée dans \cite{ja04} apporte
quelques éléments de réponse à la question \ref{min:q2}. Tout d'abord,
la réponse est négative en ce qui concerne la minoration par le
rayon d'injectivité:
\begin{theo}\label{min:tpvp}
Pour tout entier $k\geq1$ et pour toute variété $(N,h)$
telle que $b_2(N)\geq k$, il existe un fibré principal $M$ en tore $T^k$
sur $N$, une famille de métrique $(g_\varepsilon)_{\varepsilon\in]0,1]}$
sur $M$, et des réels strictement positifs $C(k,(N,h))$
et $\varepsilon_0(k,(N,h))$
tels que la courbure et le diamètre de $(M,g_\varepsilon)$ soient
uniformément bornés par rapport à $\varepsilon$,
$\Vol(M,g_\varepsilon)=\varepsilon$ pour tout $\varepsilon$, et
\begin{equation}\label{intro:eq3}
\lambda_{p,1}(M,g_\varepsilon)\leq C\cdot\injrad^{2k}(M,g_\varepsilon)
\end{equation}
pour $p=1$ et 2, et pour tout $\varepsilon<\varepsilon_0$.

 De plus, si $b_1(N)>b_2(M)$, on a aussi pour $p=2$ et 3
\begin{equation}\label{intro:eq4}
\lambda_{p,b_1(N)-b_2(M)}(M,g_\varepsilon)\leq C\cdot\injrad^{2k}
(M,g_\varepsilon)
\end{equation}
\end{theo}
 L'étude du spectre des fibrés principaux en tores qui s'effondrent, et en 
particulier
la démonstration du théorème \ref{min:tpvp}, fait appel à la théorie
des approximations diophantiennes, comme le montre l'exemple des fibrés
en tore $T^2$:

\begin{exemple}
Soit $M$ un fibré principal en tores $T^2$ sur une variété $N$ et $g$ une
métrique $T^2$-invariante sur $M$. Si on se donne un vecteur $v=(1,\alpha)
\in\R^2$ avec $\alpha$ irrationnel, il induit par l'action de $T^2$ un
champ de vecteur invariant $V$ tangent aux fibres. On définit une
famille de métrique $(g_\varepsilon)$ sur $M$ en décomposant $g$ en
la somme $g=g_V\oplus g_\bot$ d'une métrique $g_V$ dans la direction
du champ $V$
et d'une métrique $g_\bot$ sur l'orthogonal de $V$
dans $TM$, et en posant 
$g_\varepsilon=\varepsilon^2g_V\oplus g_\bot$. La variété $(M,g_\varepsilon)$ 
s'effondre sur $N$ quand $\varepsilon\to0$ car $\alpha$ est irrationnel, 
la courbure reste bornée et on montre dans \cite{ja04} qu'on peut de plus 
choisir le fibré $M$ de sorte que l'effondrement produise une petite valeur 
propre pour les $1$-formes, qui est alors de l'ordre de $\varepsilon^2$.

Si on note 
$$\mu(\alpha)=\sup\left\{\nu,\ |\alpha-\frac pq|<\frac1{q^\nu} 
\textrm{ a une infinité de solutions }(p,q)\in\Z^2\right\}$$
l'exposant d'irrationnalité de $\alpha$, on peut montrer (cf. \cite{ja05})
que 
$$\liminf_{\varepsilon\to0}\frac{\ln\lambda_{1,1}(M,g_\varepsilon)}
{\ln\injrad(M,g_\varepsilon)}=\frac{2\mu(\alpha)}{\mu(\alpha)-1}.$$ 
Rappelons que $\mu(\alpha)$ vaut 2 pour presque tous les irrationnels, et 
en particulier pour les irrationnels algébriques, mais peut aussi prendre 
toutes les valeurs dans $[2,+\infty]$ selon le choix de $\alpha$.
\end{exemple}

Le second résultat de \cite{ja04} est qu'on peut minorer le spectre
des $1$-formes d'un fibré principal en tore par le carré du volume:

\begin{theo}\label{min:th}
Soit deux réels $a$ et $d$ strictement positifs, un entier $n\geq3$ et
$(N,h)$ une variété riemannienne de dimension strictement inférieure 
à $n$. Il existe des constantes $\varepsilon_0(n,a,d,(N,h))>0$, 
$C(n,a,d,(N,h))>0$ et $C'(n,a,d,(N,h))>0$ telles que si $(M,g)$ 
est une variété riemannienne 
de dimension $n$ vérifiant $\diam(M,g)\leq d$, $|K(M,g)|\leq a$ et 
si $\pi:(M,g)\rightarrow(N,h)$ est une fibration principale de fibre 
$T^k$ qui soit une $\varepsilon$-approximation de Hausdorff avec 
$\varepsilon<\varepsilon_0$, alors
\begin{equation}\label{intro:eq2}
\lambda_{1,1}(M,g)\geq C\cdot\Vol^2(M,g)\geq C'\cdot\injrad^{2k}(M,g).
\end{equation}
\end{theo}

Le théorème \ref{min:th} soulève les deux questions suivantes qui
restent ouvertes:
\begin{question} Peut-on généraliser cette minoration aux $p$-formes
différentielles, pour tout $p$ ?
\end{question}
\begin{question}
La minoration du spectre par le volume au carré de la variété se
généralise-t-elle à d'autres familles de variétés ?
\end{question}

\section{Petites valeurs propres à courbure minorée}
Récemment a été abordé le problème de l'existence de petites valeurs
propres sous une hypothèse géométrique plus faible, à savoir 
que la courbure sectionnelle est seulement bornée inférieurement, le
diamètre restant majoré. 
Plus précisément, des exemples de petites valeurs propres ont été 
exhibés pour une famille d'effondrements à courbure minorée 
présentée par T.~Yamaguchi dans \cite{ya91}: soit $M$ une variété sur
laquelle agit un groupe de Lie compact $G$ (cette action n'est pas
nécessairement libre). On munit $M$ et $G$ des métriques bi-invariantes
$g$ et $\bar g$ respectivement, et pour tout $\varepsilon$ ont définit 
sur $M$ la métrique $g_\varepsilon$ comme étant la métrique quotient
de $((G,\varepsilon^2\bar g)\times(M,g))/G$ pour l'action diagonale de $G$.
La variété $(M,g_\varepsilon)$ tend pour la distance de Gromov-Hausdorff
vers $M/G$ quand $\varepsilon\to0$, la courbure de $g_\varepsilon$
restant uniformément minorée.

J.~Takahashi a exhibé dans \cite{ta02} un premier exemple de petite valeur 
propre
en considérant une action de $S^1$ sur $S^{2n}$ (remarque: comme la
caractéristique d'Euler des sphères de dimension paire est non nulle, 
elles ne peuvent pas s'effondrer à courbure bornée). J.~Lott
a généralisé ce résultat dans \cite{lo04} en donnant pour tous les 
effondrements de Yamaguchi décrits
ci-dessus une minoration du nombre de petites valeurs propres qui dépend
de la topologie de $M$ et $M/G$.

 Les connaissances sur ce problème restent très limitées. On ne sait par 
exemple pas si la liste de petite valeurs propres donné dans \cite{lo04} 
est exhaustive. En outre, le comportement asymptotique du spectre
dépend de la géométrie de l'effondrement, comme le montre l'exemple
des sphères de dimension impaire:

\begin{exemple}\label{cmin:ex}
On considère l'action de $\SO(2n)$ sur la sphère $S^{2n-1}$. 
L'effondrement de Yamaguchi $(S^{2n-1},g_\varepsilon)$ associé à cette
action est une simple homothétie qui effondre la sphère sur un point sans 
produire de petite valeur
propre. Cependant, comme $S^{2n-1}$ est muni d'une structure
de fibré en cercle, on peut ensuite appliquer le théorème \ref{min:cercle} 
et effondrer ce fibré de manière à obtenir une suite de métriques 
$g'_\varepsilon$ telle que $g'_\varepsilon\leq g_\varepsilon$ et 
$\lambda_{1,1}(M,g'_\varepsilon)<\varepsilon$, la courbure restant 
uniformément minorée.

\end{exemple}

\begin{remarque}
L'exemple \ref{cmin:ex} peut se généraliser aux fibrés en cercles
dont la base s'effondre à courbure minorée. On peut noter par ailleurs
qu'on utilise le fait que la variété admet une petite valeur propre
à courbure bornée.
\end{remarque}

\begin{question}
Existe-t-il une variété $M$ dont le volume minimal est nul, qui 
n'admet pas de petite valeur propre à courbure sectionnelle bornée mais 
qui en admet à courbure minorée?
\end{question}

 On peut aussi envisager d'affaiblir encore l'hypothèse sur la courbure:
\begin{question}
Existe-t-il une variété qui n'admet pas de petite valeur propre à courbure
sectionnelle minorée mais qui en admet à courbure de Ricci minorée?
\end{question}

 Enfin, on peut reformuler la question \ref{min:q2} avec ces hypothèses:
\begin{question}
À diamètre borné et courbure minorée, peut-on minorer la première
valeur propre du laplacien de Hodge-de~Rham par le volume de la
variété au carré ?
\end{question}

\noindent Pierre \textsc{Jammes}\\
Université d'Avignon\\
laboratoire de mathématiques\\
33 rue Louis Pasteur\\
F-84000 Avignon\\
\texttt{Pierre.Jammes@univ-avignon.fr}
\end{document}